\documentclass[11pt]{article}
\usepackage{enumerate}
\usepackage{amssymb,a4wide,latexsym,makeidx,epsfig,fleqn}
\usepackage{amsthm}
\usepackage{amsmath}
\usepackage{enumerate}
\newtheorem{theorem}{Theorem}[section]

\newtheorem{lemma}[theorem]{Lemma}

\begin{document}
\textwidth 150mm \textheight 225mm
\title{Minimal Skew energy of oriented bicyclic graphs with a given diameter
\thanks{ Supported by
the National Natural Science Foundation of China (Nos.11171273 and 11601431) and the Seed Foundation of Innovation and Creation for Graduate Students in Northwestern Polytechnical University (No.Z2016170).}}
\author{{Xiangxiang Liu and Ligong Wang\footnote {Corresponding author.}
 }\\
{\small Department of Applied Mathematics, School of Science, Northwestern
Polytechnical University,}\\ {\small  Xi'an, Shaanxi 710072,
People's Republic
of China. }\\{\small E-mail: xxliumath@163.com, lgwangmath@163.com}}
\date{}
\maketitle
\begin{center}
\begin{minipage}{120mm}
\vskip 0.3cm
\begin{center}
{\small {\bf Abstract}}
\end{center}
{\small Let $S(G^{\sigma})$ be the skew-adjacency matrix of the oriented graph $G^{\sigma}$, which is obtained from a simple undirected graph $G$ by assigning an orientation $\sigma$ to each of its edges. The skew energy of an oriented graph $G^{\sigma}$ is defined as the sum of absolute values of all eigenvalues of $S(G^{\sigma})$. For any positive integer $d$ with $3\leq d\leq n-3$, we determine the graph with minimal skew energy among all oriented bicyclic graphs that contain no vertex disjoint odd cycle of lengths $s$ and $l$ with $s+l\equiv 2(mod 4)$ on $n$ vertices with a given diameter $d$.

 }
\end{minipage}
\end{center}

\section{Introduction }
\label{sec:ch6-introduction}

Let $G$ be a simple undirected graph with an orientation $\sigma$, which assigns to each edge a direction so that $G^{\sigma}$ becomes an oriented graph. Then $G$ is usually called the underlying graph of $G^{\sigma}$. The skew-adjacency matrix of $G^{\sigma}$ with vertex set $V(G)=\{1,2,...,n\}$ is the $n\times n$ matrix $S(G^{\sigma})=[s_{ij}],$ $s_{ij}=1$ and $s_{ji}=-1$ if $(i,j)$ is an arc of $G^{\sigma}$, and $s_{ij}=s_{ji}=0$ otherwise. Since $S(G^{\sigma})$ is a real skew symmetric matrix, all eigenvalues $\{\lambda_{1},\lambda_{2},...,\lambda_{n}\}$ of $S(G^{\sigma})$ are pure imaginary numbers or $0$.

The skew energy of an oriented graph $G^{\sigma}$, denoted by $E_{s}(G^{\sigma})$, is defined as the sum of absolute values of all eigenvalues of $S(G^{\sigma})$ (see \cite{ABS}), that is

$$E_{s}(G^{\sigma})=\sum\limits^{n}_{i=1}|\lambda_{i}|,$$\\where $\lambda_{1},\lambda_{2},...,\lambda_{n}$ are the eigenvalues of the skew-adjacency matrix $S(G^{\sigma})$, namely the $n$ roots of $\phi(G^{\sigma};x)=0$. Here $\phi(G^{\sigma};x)=det(xI_{n}-S(G^{\sigma}))=\sum^{n}_{i=0}a_{i}(G^{\sigma})x^{n-i}$ is the skew characteristic polynomial of $G^{\sigma}$, where $I_{n}$ is the unit matrix of order $n$. Since $S(G^{\sigma})$ is a real skew symmetric matrix, $a_{2i}(G^{\sigma})\geq0$ and $a_{2i+1}(G^{\sigma})=0$ for all $0\leq i\leq\lfloor\frac{n}{2}\rfloor$ (see \cite{ABS}). So we have

$$\phi(G^{\sigma};x)=\sum\limits^{\lfloor\frac{n}{2}\rfloor}_{i=0}a_{2i}(G^{\sigma})x^{n-2i}.\eqno{(1)}$$\\By using the coefficients of $\phi(G^{\sigma};x)$, the skew energy $E_{s}(G^{\sigma})$ can be expressed by the following integral formula \cite{HSZ}

$$E_{s}(G^{\sigma})=\frac{2}{\pi}\int^{+\infty}_{0}\frac{1}{x^{2}}\ln\left[\sum\limits^{\lfloor\frac{n}{2}\rfloor}_{i=0}a_{2i}(G^{\sigma})x^{2i}\right]dx. \eqno{(2)} $$
\\It follows that $E_{s}(G^{\sigma})$ is a strictly monotonously increasing function of $a_{2i}(G^{\sigma})$ for $0\leq i\leq\lfloor\frac{n}{2}\rfloor$ for any oriented graph. Note that $a_{0}(G^{\sigma})=1$ and $a_{2}(G^{\sigma})$ equals to the number of the edges in $G$. This provides a useful way for comparing the skew energies of a pair of oriented graphs.

Let $G_{1}^{\sigma_{1}}$ and $G_{2}^{\sigma_{2}}$ be two oriented graphs of order $n$. If $a_{2i}(G_{1}^{\sigma_{1}})\leq a_{2i}(G_{2}^{\sigma_{2}})$ for all $0\leq i\leq\lfloor\frac{n}{2}\rfloor$, we write $G_{1}^{\sigma_{1}}\preceq G_{2}^{\sigma_{2}}.$
Furthermore, if $G_{1}^{\sigma_{1}}\preceq G_{2}^{\sigma_{2}}$ and there exists at least one index $j$ such that $a_{2j}(G_{1}^{\sigma_{1}})<a_{2j}(G_{2}^{\sigma_{2}})$, we write $G_{1}^{\sigma_{1}}\prec G_{2}^{\sigma_{2}}.$ If $a_{2i}(G_{1}^{\sigma_{1}})=a_{2i}(G_{2}^{\sigma_{2}})$ for all $0\leq i\leq\lfloor\frac{n}{2}\rfloor$, we write $G_{1}^{\sigma_{1}}\sim G_{2}^{\sigma_{2}}.$ According to the integral formula (2), we have

$G_{1}^{\sigma_{1}}\preceq G_{2}^{\sigma_{2}}\Rightarrow E_{s}(G_{1}^{\sigma_{1}})\leq E_{s}(G_{2}^{\sigma_{2}})$;

$G_{1}^{\sigma_{1}}\prec G_{2}^{\sigma_{2}}\Rightarrow E_{s}(G_{1}^{\sigma_{1}})< E_{s}(G_{2}^{\sigma_{2}})$.

The study on the extremal values of energy for oriented graphs is of importance for the chemical graph theory, and a lot of interesting results have been reported. For the oriented unicyclic graphs of order $n$, Hou et al. \cite{HSZ} obtained the oriented graphs with the 1st-minimal, the 2nd-minimal and the maximal skew energies, and Zhu \cite{ZJM} determined the oriented graphs with the first $\lfloor\frac{n-9}{2}\rfloor$ largest skew energies. For the oriented bicyclic graphs, Shen et al. \cite{SHZ} deduced the oriented graphs with the minimal and maximal skew energies, and Wang et al. \cite{WZY} characterized the oriented graph with the second largest skew energy. Zhu and Yang \cite{ZYJ} obtained the oriented unicyclic graphs that have perfect matchings with the minimal skew energy. Yang et al. \cite{YGX} determined the oriented unicyclic graphs of a fixed diameter with the minimal skew energy. Some other results about the extremal skew energies can be found in Refs. \cite{TGX,CLL,GXS}. For a survey on skew energy of oriented graphs, one can refer to \cite{LLX}.

This paper is organized as follows: In Section 2, we give some notations and preliminary results, which will be used in the following discussion. The graph with minimal skew energy among  all oriented bicyclic graphs that contain no vertex disjoint odd cycle of lengths $s$ and $l$ with $s+l\equiv 2(mod 4)$ on $n$ vertices with a given diameter $d$ will be determined in Section 3, where $3\leq d\leq n-3$. For $d=2$, we can refer to \cite{SHZ}.

\section{Preliminary Results}
\label{sec:ch-sufficient}

Let $G=(V(G),E(G))$ be a simple graph. Denote by $G-e$ the graph obtained from $G$ by deleting the edge $e$ and by $G-v$ the graph obtained from $G$ by deleting the vertex $v$ together with all edges incident to it. Let $d(G)$ be the diameter of $G$, which is defined as the greatest distance between any two vertices in $G$. The union of the graphs $G_{1}=(V(G_{1}),E(G_{1}))$ and $G_{2}=(V(G_{2}),E(G_{2}))$, denoted by $G_{1}\cup G_{2}$, is the graph with vertex set $V(G_{1})\cup V(G_{2})$ and edge set $E(G_{1})\cup E(G_{2})$. $N(u)$ denotes the neighborhood of $u$. We refer to Cvetkovi\'{c} et al. \cite{DMH} for undefined terminology and notation.

For convenience, in terms of defining subgraph, matching, degree, diameter, etc., of an oriented graph, we focus only on its underlying graph. Moreover, we will briefly use the notations $S_{n}$, $P_{n}$ and $C_{n}$ to denote the oriented star, the oriented path and the oriented cycle on $n$ vertices, respectively, if no conflict exists there.

Let $C$ be an even cycle of $G$. Then we say $C$ is evenly oriented relative to $G^{\sigma}$ if it has even number of edges oriented in the
direction of the routing, otherwise $C$ is oddly oriented. A linear subgraph $L$ of $G$ is a disjoint union of some edges and some cycles in $G$. A linear subgraph $L$ is called evenly linear subgraph if the number of vertices of $L$ is even. $\varepsilon\mathcal{L}_{i}$  denotes the set of all evenly linear subgraph of $G$ with $i$ vertices.

Let $G^{\sigma}$ be an oriented graph of $G$. Let $W$ be a subset of $V(G)$ and $\overline{W}=V(G)\setminus W$. The orientation $G^{\tau}$ of $G$ obtained from $G^{\sigma}$ by reversing the orientations of all arcs between $W$ and $\overline{W}$. Then $G^{\tau}$ is said to be obtained from $G^{\sigma}$ by a switching with respect to $W$. Moreover, two orientations $G^{\sigma}$ and $G^{\tau}$ of a graph $G$ are said to be switching equivalent if $G^{\tau}$ can be obtained from $G^{\sigma}$ by a sequence of switchings. As noted in \cite{ABS}, since the skew-adjacency matrices obtained by a switching are similar, their skew energies are equal.

It is easy to verify that up to switching equivalence there are just two orientations of a cycle $C$: $(1)$ Just one edge on the cycle has the opposite orientation to that of others, we denote this orientation by $+$. $(2)$ All edges on the cycle $C$ have the same orientation, we denote this orientation by $-$. So if a cycle is of even length and oddly oriented, then it is equivalent to the orientation $+$. If a cycle is of even length and evenly oriented, then it is equivalent to the orientation $-$.

 Adiga et al. \cite{ABS} showed that the skew energy of a directed tree is independent of its orientation, which is equal to the energy of its underlying tree. So by switching equivalence, for a unicyclic oriented graph or bicyclic oriented graph, we only need to consider the orientations of cycles.

Let $C_{a}$, $C_{b}$ be two cycles in bicyclic graph $G$ with $t(t\geq 0)$ common vertices. If $t\leq1$, then $G$ contains exactly two cycles. If $t\geq2$, then $G$ contains exactly three cycles. The third cycle is denoted by $C_{c}$, where $c=a+b-2t+2$. Let $C_{a}=v_{0}v_{1}\cdots v_{a-1}v_{0}$ and $C_{b}=u_{0}u_{1}\cdots u_{b-1}u_{0}$. If $C_{a}$ and $C_{b}$ have no common vertices, then $C_{a}$ and $C_{b}$ are connected by a path $P$, say from $v_{0}$ to $u_{0}$. Let $l(G)$ be the length of $P$. If $t\geq 1$, $C_{c}=u_{0}u_{b-1}\cdots u_{t}u_{t-1}v_{t}v_{t+1}\cdots v_{a-1}v_{0}$ is the third cycle, where $v_{0}=u_{0}$, $v_{1}=u_{1}$,$\cdots$ $v_{t-1}=u_{t-1}$. If we write $w_{0}=u_{0},w_{1}=u_{b-1},\cdots w_{c-1}=v_{a-1}$, then $C_{c}=w_{0}w_{1}\cdots w_{c-1}w_{0}$.

For convenience, we denote by $G^{+}$ (resp. $G^{-}$) the unicyclic graph on which the orientation of a cycle is of orientation $+$ (resp. $-$), and denote by $G^{\ast}$ the unicyclic graph on which the orientation of a cycle is of arbitrary orientation $\ast$. If $t\leq1$, we denote by $G^{\alpha,\beta}$ the bicyclic graph on which cycle $C_{a}$ is of orientation $\alpha$ and cycle $C_{b}$ is of orientation $\beta$,  where $\alpha,\beta\in\{+,-,\ast\}$. If $t\geq 2$, we denote by $G^{\alpha,\beta,\gamma}$ the bicyclic graph on which $C_{a}$ is of orientation $\alpha$, $C_{b}$ is of orientation $\beta$ and $C_{c}$ is of orientation $\gamma$, where $\alpha,\beta,\gamma\in\{+,-,\ast\}$.

The following results are the cornerstone of our discussion below, which gives an interpretation of all coefficients of the skew characteristic polynomial of an oriented graph.

\noindent\begin{lemma}\label{le:2.1} (\cite{HLY})
 Let $G^{\sigma}$ be an oriented graph of a graph $G$ with the skew characteristic polynomial $\phi(G^{\sigma};x)=\sum\limits^{n}_{i=0}a_{i}(G^{\sigma})x^{n-i}.$ Then \\$$a_{i}(G^{\sigma})=\sum\limits_{L\in \varepsilon\mathcal{L}_{i}}(-2)^{p_{e}(L)}2^{p_{o}(L)},$$
\\where $p_{e}(L)\ (resp.\ p_{o}(L))$ is the number of all evenly (resp. oddly) oriented cycles of a linear subgraph $L$ relative to $G^{\sigma}$.
\end{lemma}

\noindent\begin{lemma}\label{le:2.2} (\cite{HLY})
 Let $e=(u,v)$ be an arc of an oriented graph $G^{\sigma}$. Then
 \begin{align*}
 a_{i}(G^{\sigma})=&a_{i}(G^{\sigma}-e)+a_{i-2}(G^{\sigma}-u-v)+2\sum\limits_{e\in C\in Od(G^{\sigma})}a_{i-|V(C)|}(G^{\sigma}-V(C))\\
                     &-2\sum\limits_{e\in C\in Ev(G^{\sigma})}a_{i-|V(C)|}(G^{\sigma}-V(C)),
\end{align*}
where $Od(G^{\sigma})\ (resp.\ Ev(G^{\sigma}))$ denotes the set of all oddly (resp. evenly) cycles of $G^{\sigma}$.
\end{lemma}

\noindent\begin{lemma}\label{le:2.3} (\cite{XGH})
 Let $v$ be a vertex of an oriented graph $G^{\sigma}$. Then
 \begin{align*}
 a_{i}(G^{\sigma})=&a_{i}(G^{\sigma}-v)+\sum\limits_{u\in N(v)}a_{i-2}(G^{\sigma}-u-v)+2\sum\limits_{v\in C\in Od(G^{\sigma})}a_{i-|V(C)|}(G^{\sigma}-V(C))\\
                     &-2\sum\limits_{v\in C\in Ev(G^{\sigma})}a_{i-|V(C)|}(G^{\sigma}-V(C)),
\end{align*}
where $Od(G^{\sigma})\ (resp.\ Ev(G^{\sigma}))$ denotes the set of all oddly (resp. evenly) cycles of $G^{\sigma}$.
\end{lemma}

From the Lemma 2.2, we can obtain easily Lemmas 2.4 and 2.5.

\noindent\begin{lemma}\label{le:2.4}
Let $e$ be a cut edge of $G$. Then $G^{\sigma}\succeq G^{\sigma}-e$
\end{lemma}

\noindent\begin{lemma}\label{le:2.5}
 Let $G$ be a unicyclic graph with $n$ vertices or a bicyclic graph with $n$ vertices. Then
$G^{\sigma}\succeq S_{n}$.
\end{lemma}

Let $\mathcal{T}(n,d)$ be the class of trees with $n$ vertices and diameter $d$. Denote by $T_{n,d}$ the tree obtained from the path $P_{d-1}$ and the star $S_{n-d+2}$ by identifying one pendent vertex of them. $\mathcal{U}(n,d)$ denotes the class of unicyclic graphs with $n$ vertices and diameter $d$ and $U_{n,d}$ denotes the unicyclic graph obtained from the cycle $C_{4}$ by attaching a pendent vertex of the path $P_{d-2}$ and $n-d-1$ pendent edges to its two non-adjacent vertices respectively; see Figure 1.

Let $\mathcal{B}(n)$ be the class of bicyclic graphs with $n$ vertices and contains no vertex disjoint odd cycles of lengths $s$ and $l$ with $s+l\equiv 2(mod 4)$. Let $\mathcal{B}(n,d)$ be the class of bicyclic graphs $\mathcal{B}(n)$ with diameter $d$ where $2\leq d\leq n-3$. Denote by $B_{n,d}$ the bicyclic graph obtained from the $K_{2,3}$ by attaching a pendent vertex of the path $P_{d-2}$ and $n-d-2$ pendent edges to its two vertices of degree three respectively; see Figure 1.

We have known that the skew energy of a directed tree is independent of its orientation. Hence, the following results for undirected trees apply equally well to oriented trees, which will be cited in the following discussion directly.

\noindent\begin{lemma}\label{le:2.6} (\cite{GPI})
For $n\geq2,$ $P_{n}\succeq P_{i}\cup P_{n-i}\succeq P_{1}\cup P_{n-1}$.
\end{lemma}

\noindent\begin{lemma}\label{le:2.7} (\cite{GIU})
Let $n\geq5,$ $T_{n}$ denote any tree with order $n$ and $T_{n}\neq P_{n},S_{n}$. Then $P_{n}\succeq T_{n}\succeq S_{n}$.
\end{lemma}

\noindent\begin{lemma}\label{le:2.8} (\cite{YYW})
Let $T\in \mathcal{T}(n,d)$. Then $T\succeq T_{n,d}$.
\end{lemma}

\noindent\begin{lemma}\label{le:2.9} (\cite{LZF})
If $d>d_{0}\geq3,$ then $T_{n,d}\succeq T_{n,d_{0}}$.
\end{lemma}

\noindent\begin{lemma}\label{le:2.10} (\cite{YZY})
If $2\leq d_{1}< n_{1}-2,$ then $T_{n_{1},d_{1}}\cup T\succeq T_{n_{1}+n_{2}-1,d_{1}+d_{2}}$, where $T=T_{n_{2},d_{2}}$ if $2\leq d_{2}< n_{2}-1$ or $P_{2}$ if $n_{2}=2$ and $d_{2}=1$.
\end{lemma}

\noindent\begin{lemma}\label{le:2.11} (\cite{YGX})
Let $U\in \mathcal{U}(n,d)$ with  $n\geq6$ and $3\leq d\leq n-2$. Then $U^{\sigma}\succeq U^{-}_{n,d}$.
\end{lemma}

\noindent\begin{lemma}\label{le:2.12}
For $3\leq d\leq n-2$,  $U^{-}_{n,d}\succeq T_{n,d}$.
\end{lemma}

{\bf Proof.} By Lemma 2.2,
\begin{align*}
 a_{2i}(U^{-}_{n,d})=&a_{2i}(T_{n,d})+a_{2i-2}(P_{d-3}\cup S_{n-d+1})-2a_{2i-4}(P_{d-3})\\
                     =&a_{2i}(T_{n,d})+a_{2i-2}(P_{d-3}\cup S_{n-d-1})\geq a_{2i}(T_{n,d}).
\end{align*}
\hfill$\blacksquare$
\begin{center}
\includegraphics [width=10 cm, height=4 cm]{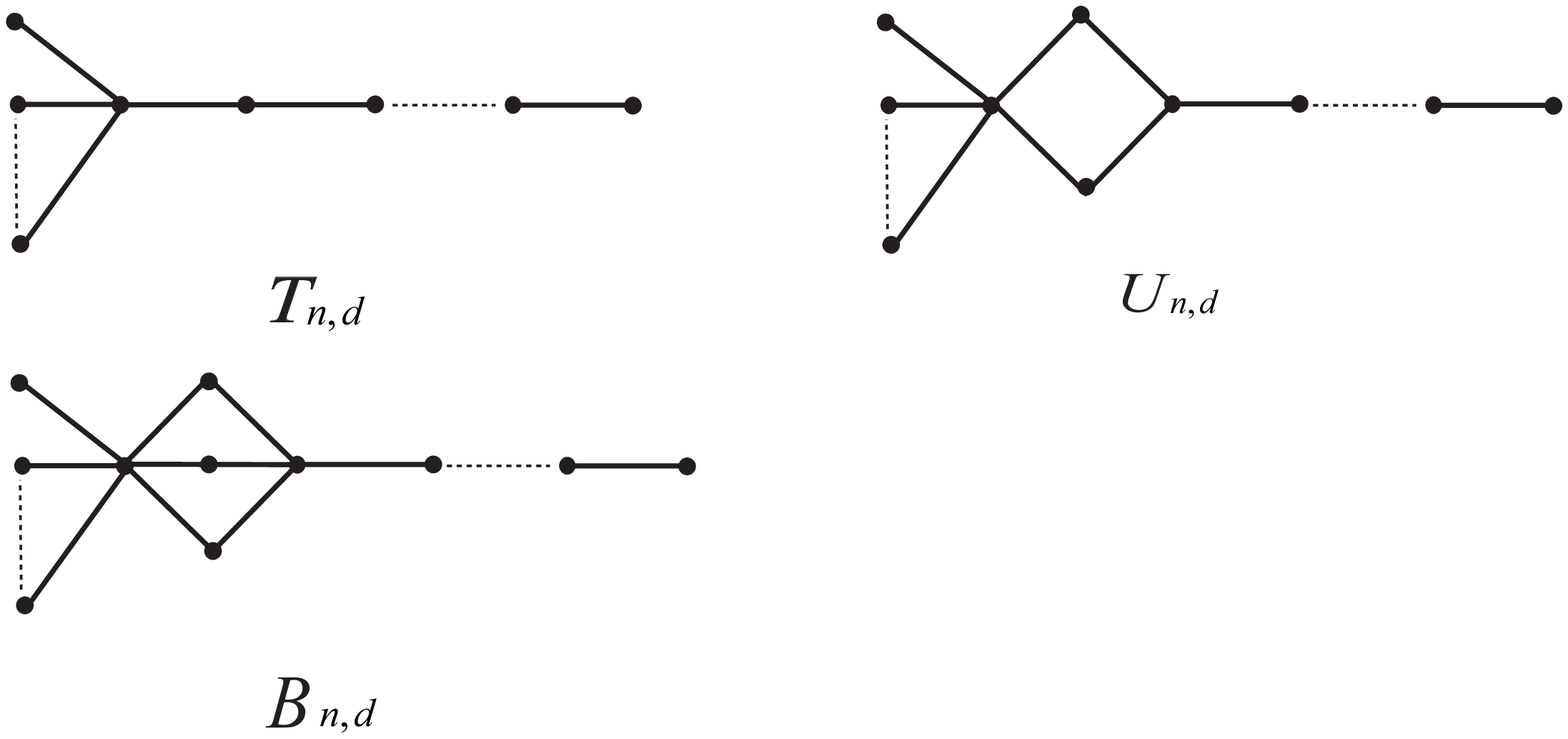}
\centerline{Figure 1: The tree $T_{n,d}$, the unicyclic graph $U_{n,d}$ and the bicyclic graph $B_{n,d}$.}
\end{center}

\noindent\begin{lemma}\label{le:2.13}
If $3\leq d_{0} < d\leq n-2$, then $U^{-}_{n,d}\succeq U^{-}_{n,d_{0}}$.
\end{lemma}

{\bf Proof.} If $d=4$, we have $U^{-}_{n,4}\succeq U^{-}_{n,3}$ by Lemma 2.1. If $d\geq5$, by Lemmas 2.2, 2.4 and 2.12,
\begin{align*}
 a_{2i}(U^{-}_{n,d})=&a_{2i}(U^{-}_{n-1,d-1})+a_{2i-2}(U^{-}_{n-2,d-2})\\
                     \geq&a_{2i}(U^{-}_{n-1,d-1})+a_{2i-2}(T_{n-2,d-2})\\
                     \geq&a_{2i}(U^{-}_{n-1,d-1})+a_{2i-2}(T_{d-1,d-3})\\
                     =&a_{2i}(U^{-}_{n,d-1}),
\end{align*}
so $U^{-}_{n,d}\succeq U^{-}_{n,d-1}\succeq \cdots \succeq U^{-}_{n,d_{0}}$. \hfill$\blacksquare$

Similarly, we have the following result.

\noindent\begin{lemma}\label{le:2.14}
If $3\leq d_{0} < d\leq n-2$, then $B^{-,-,-}_{n,d}\succeq B^{-,-,-}_{n,d_{0}}$.
\end{lemma}

\section{Main Results}
\label{sec:ch-inco}

\noindent\begin{lemma}\label{le:3.1}
Let $G\in \mathcal{B}(n,n-3)$ with $n\geq 6$, and $G\neq B_{n,n-3}$. Then $G^{\sigma}\succ B^{-,-,-}_{n,n-3}$.
\end{lemma}

{\bf Proof.} We prove this lemma by induction on $n$.

When $n=6$, then $G\in \mathcal{B}(6,3)$ and $G\neq B_{6,3}$. Then $G$ is isomorphic to one of the graphs in Figure 2. By Lemma 2.1, we have $G^{\sigma}\succ B^{-,-,-}_{6,3}$.

\begin{center}
\includegraphics [width=11 cm, height=5 cm]{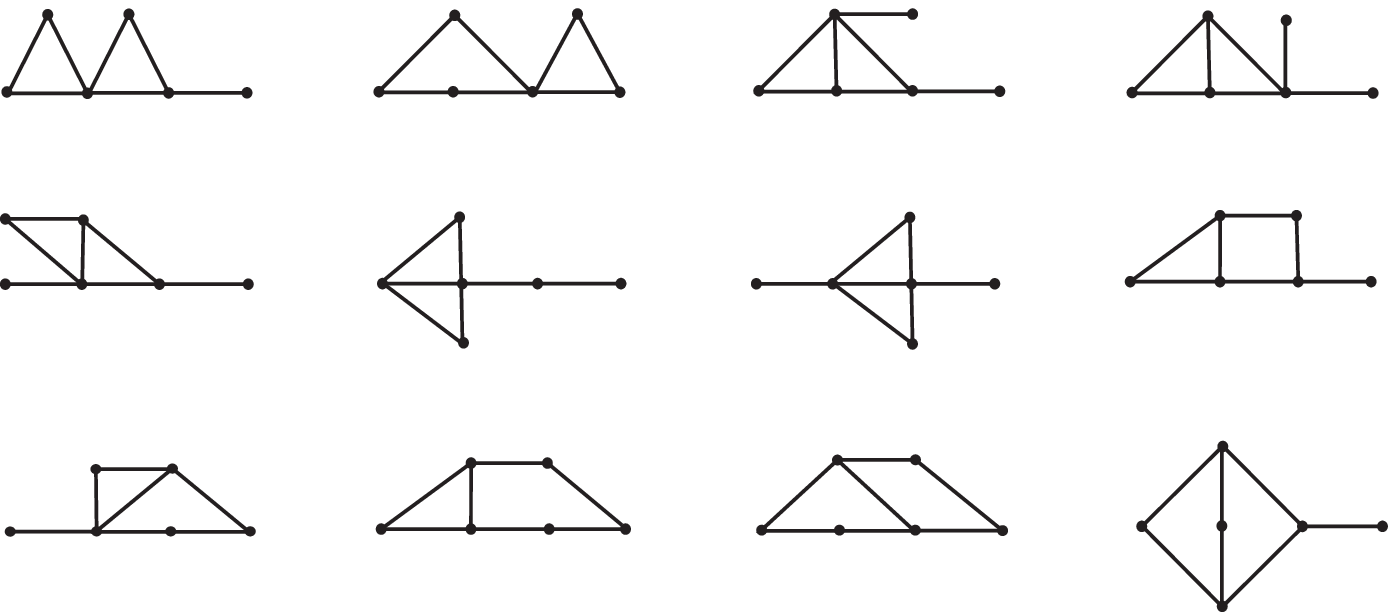}
\centerline{Figure 2: Graphs in $\mathcal{B}(6,3)$ except $B_{6,3}$.}
\end{center}

When $n=7$, then $G\in \mathcal{B}(7,4)$ and $G\neq B_{7,4}$. Then $G$ is isomorphic to one of the graphs in Figure 3. By Lemma 2.1, we have $\phi(B^{-,-,-}_{7,4};x)=x^{7}+8x^{5}+7x^{3}$. By a directly calculation, we have $a_{4}(G^{\sigma})>a_{4}(B^{-,-,-}_{7,4})=7$. So
$G^{\sigma}\succ B^{-,-,-}_{7,4}$.

Suppose that the result holds for graphs of $\mathcal{B}(n-1,n-4)$ and $\mathcal{B}(n-2,n-5)$ with $n\geq8$. Now suppose that $G\in \mathcal{B}(n,n-3)$ and $G\neq B_{n,n-3}$.

\begin{center}
\includegraphics [width=15 cm, height=8 cm]{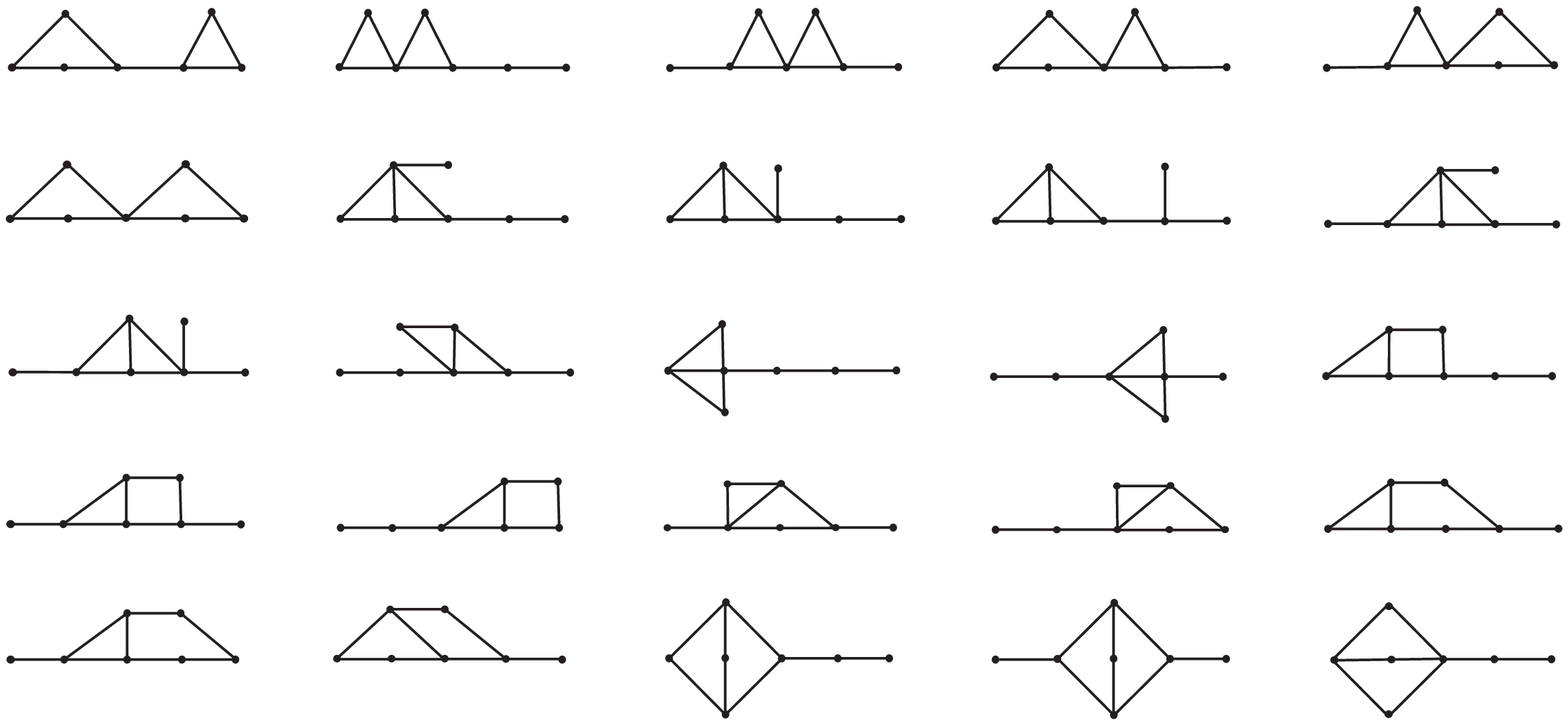}
\centerline{Figure 3: Graphs in $\mathcal{B}(7,4)$ except $B_{7,4}$.}
\end{center}

{\bf Case 1.} There exists a pendent vertex $u$ in $G$ such that the degree of its neighbor
$v$ is two. Then $G-u\in \mathcal{B}(n-1,n-4)$ and $G-u-v\in \mathcal{B}(n-2,n-5)$. By Lemma 2.2, we have

$a_{2i}(G^{\sigma})=a_{2i}(G^{\sigma}-u)+a_{2i-2}(G^{\sigma}-u-v)$,

$a_{2i}(B^{-,-,-}_{n,n-3})=a_{2i}(B^{-,-,-}_{n-1,n-4})+a_{2i-2}(B^{-,-,-}_{n-2,n-5})$.
\\Note that $G\neq B_{n,n-3}$, thus $G-u\neq B_{n-1,n-4}$ and $G-u-v\neq B_{n-2,n-5}$. Combining with the induction hypothesis, then $G^{\sigma}\succ B^{-,-,-}_{n,n-3}$.

\begin{center}
\includegraphics [width=15 cm, height=2 cm]{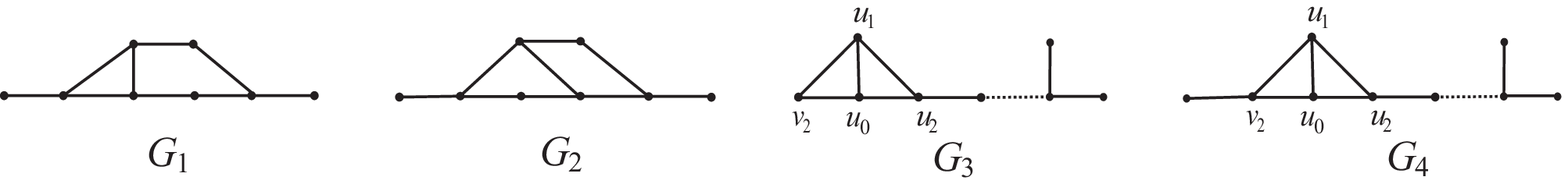}
\centerline{Figure 4: Graphs $G_{j}$ $j$=1,2,3,4.}
\end{center}

{\bf Case 2.} The neighbor of any pendent vertex has degree at least three or there
is no pendent vertex. Then $G$ is isomorphic to some $G_{j}$ in Figure 4, $j=1,2,3,4$, or $G$ contains one quadrangle which has at most one common vertex with another cycle that is a triangle or a quadrangle. For $n=8$, if $G$ is isomorphic to $G_{1}$ or $G_{2}$, then $G^{\sigma}\succ B^{-,-,-}_{8,5}$ by a directly calculation.

{\bf Subcase 2.1.} $G$ is isomorphic to $G_{3}$ or $G_{4}$, then by Lemmas 2.2 and 2.8, we have
\begin{align*}
 a_{2i}(G^{\sigma})=&a_{2i}(G^{\sigma}-u_{0}u_{2})+a_{2i-2}(G^{\sigma}-u_{0}-u_{2})-2a_{2i-4}(G^{\sigma}-u_{0}-u_{1}-u_{2}-v_{2})\\
                     =&a_{2i}(G^{\sigma}-u_{0}u_{2}-u_{0}v_{2})+a_{2i-2}(G^{\sigma}-u_{0}-v_{2})+a_{2i-2}(G^{\sigma}-u_{0}-u_{2})\\
                     &-2a_{2i-4}(G^{\sigma}-u_{0}-u_{1}-u_{2}-v_{2})\\
                     =&a_{2i}(G^{\sigma}-u_{0}u_{2}-u_{0}v_{2})+a_{2i-2}(G^{\sigma}-u_{0}-v_{2}-u_{1}u_{2})\\
                     &+a_{2i-2}(G^{\sigma}-u_{0}-u_{2}-u_{1}v_{2})\\
                     \geq&a_{2i}(T_{n,n-3})+2a_{2i-2}(P_{n-6})=a_{2i}(B^{-,-,-}_{n,n-3}),
\end{align*}
thus $G^{\sigma}\succ B^{-,-,-}_{n,n-3}$.

{\bf Subcase 2.2.} $G$ contains one quadrangle which has at most one common vertex with another cycle that is a triangle or a quadrangle. If $n=8, 9$, it can be checked by Lemma 2.1 that $G^{\sigma}\succ B^{-,-,-}_{n,n-3}$. If $n\geq10$, suppose that $C_{b}=u_{0}u_{1}u_{2}u_{3}u_{0}$ (see Section 2) is a quadrangle. By Lemmas 2.2, 2.11 and 2.12, we have
\begin{align*}
 a_{2i}(G^{\sigma})=&a_{2i}(G^{\sigma}-u_{0}u_{1})+a_{2i-2}(G^{\sigma}-u_{0}-u_{1})-2a_{2i-4}(G^{\sigma}-u_{0}-u_{1}-u_{2}-u_{3})\\
                     =&a_{2i}(G^{\sigma}-u_{0}u_{1}-u_{1}u_{2})+a_{2i-2}(G^{\sigma}-u_{1}-u_{2})+a_{2i-2}(G^{\sigma}-u_{0}-u_{1})\\
                     &-2a_{2i-4}(G^{\sigma}-u_{0}-u_{1}-u_{2}-u_{3})\\
                     =&a_{2i}(G^{\sigma}-u_{0}u_{1}-u_{1}u_{2})+a_{2i-2}(G^{\sigma}-u_{1}-u_{2}-u_{0}u_{3})\\
                     &+a_{2i-2}(G^{\sigma}-u_{0}-u_{1}-u_{2}u_{3})\\
                     \geq&a_{2i}(U^{-}_{n-1,n-3})+a_{2i-2}(U^{-}_{n-3,n-5})+a_{2i-2}(U^{-}_{n-4,n-6})\\
                     \geq&a_{2i}(U^{-}_{n-1,n-3})+a_{2i-2}(U^{-}_{n-5,n-7}\cup P_{2})+a_{2i-2}(U^{-}_{n-4,n-6}).
\end{align*}
By Lemma 2.2, we have $a_{2i}(B^{-,-,-}_{n,n-3})=a_{2i}(U^{-}_{n-1,n-3})+a_{2i-2}(P_{n-6}\cup P_{2})+a_{2i-2}(P_{n-5})$, so $G^{\sigma}\succ B^{-,-,-}_{n,n-3}$.
 \hfill$\blacksquare$

\noindent\begin{lemma}\label{le:3.2}
Let $G\in \mathcal{B}(n,d)$ with $3\leq d\leq n-4$. If $G$ contains no pendent vertices, then $G^{\sigma}\succ B^{-,-,-}_{n,d}$.
\end{lemma}

{\bf Proof.} Let $b\geq a$. Since $d\leq n-4$, we have $b\geq 5$. By Lemma 2.2, we have
$$a_{2i}(B^{-,-,-}_{n,d+1})=a_{2i}(U^{-}_{n,d+1})+a_{2i-2}(P_{d-2}\cup S_{n-d})-4a_{2i-4}(P_{d-2}).\eqno{(3)}$$
$$a_{2i}(B^{-,-,-}_{n,d+1})=a_{2i}(U^{-}_{n-1,d+1})+a_{2i-2}(P_{d-2}\cup S_{n-d-2})+a_{2i-2}(P_{d-1}).  \quad (By\ (3))\eqno{(4)}$$
$$a_{2i}(B^{-,-,-}_{n,d+1})=a_{2i}(T_{n-1,d+1})+a_{2i-2}(P_{d})+2a_{2i-2}(P_{d-2}\cup S_{n-d-3}). \quad  (By\ (4))\eqno{(5)}$$
$$a_{2i}(B^{-,-,-}_{n,d+1})=a_{2i}(T_{n,d+1})+2a_{2i-2}(P_{d-2}\cup S_{n-d-3}). \quad  (By\ (3))\eqno{(6)}$$
{\bf Case 1.} When $t\leq 1$, there are exactly two cycles $C_{a}$ and $C_{b}$ in $G$. Then $d=\lfloor \frac {a}{2}\rfloor+\lfloor \frac {b}{2}\rfloor+l(G)$ (see Section 2 for $l(G)$).

{\bf Subcase 1.1.} The length of $C_{b}$ is odd. Then
$d(G-u_{1}u_{2})=\lfloor \frac {a}{2}\rfloor+b+l(G)-2\geq d+1$, $d(G-u_{1}-u_{2})\geq d$ and $d\geq 3$. By Lemmas 2.2, 2.11, 2.12 and 2.13, we have
\begin{align*}
 a_{2i}(G^{\sigma})=&a_{2i}(G^{\sigma}-u_{1}u_{2})+a_{2i-2}(G^{\sigma}-u_{1}-u_{2})\\
                     \geq&a_{2i}(U^{-}_{n,d+1})+a_{2i-2}(U^{-}_{n-2,d})\\
                     \geq&a_{2i}(U^{-}_{n,d+1})+a_{2i-2}(T^{-}_{n-2,d}).
\end{align*}
Combining with $(3)$ and Lemma 2.14, then $G^{\sigma}\succ B^{-,-,-}_{n,d+1}\succ B^{-,-,-}_{n,d}$.

{\bf Subcase 1.2.}  The length of $C_{b}$ is even. Then $b\geq 6$. Hence $d(G-u_{1}u_{2}-u_{2}u_{3})=\lfloor \frac {a}{2}\rfloor+b+l(G)-3\geq d$, $d(G-u_{2}-u_{3}-u_{4}u_{5})\geq d-2$, $d(G-u_{1}-u_{2}-u_{3}u_{4})\geq d-1$. and $d\geq4$. If $d=4$, it can be checked by Lemmas 2.1 that $G^{\sigma}\succ B^{-,-,-}_{n,4}$. If $d\geq5$, by Lemmas 2.2, 2.11, 2.12 and 2.13, we have
\begin{align*}
 a_{2i}(G^{\sigma})\geq&a_{2i}(G^{\sigma}-u_{1}u_{2}-u_{2}u_{3})+a_{2i-2}(G^{\sigma}-u_{1}-u_{2})\\
                     &+a_{2i-2}(G^{\sigma}-u_{2}-u_{3})-2a_{2i-b}(G^{\sigma}-V(C_{b}))\\
                     \geq&a_{2i}(G^{\sigma}-u_{1}u_{2}-u_{2}u_{3})+a_{2i-2}(G^{\sigma}-u_{1}-u_{2}-u_{3}u_{4})\\
                     &+a_{2i-2}(G^{\sigma}-u_{2}-u_{3}-u_{4}u_{5})\\
                     \geq&a_{2i}(U^{-}_{n-1,d})+a_{2i-2}(U^{-}_{n-3,d-1})+a_{2i-2}(U^{-}_{n-3,d-2})\\
                      \geq&a_{2i}(U^{-}_{n-1,d})+a_{2i-2}(T_{n-3,d-1})+a_{2i-2}(T_{n-3,d-2}).
\end{align*}
Combining with $(4)$, $G^{\sigma}\succ B^{-,-,-}_{n,d}$.

{\bf Case 2.} $t\geq 2$. 
Note that $a-t+1\geq t-1$ and $b-t+1\geq t-1$. Then $c\geq b$ and $d=\lfloor \frac{c}{2}\rfloor=\lfloor \frac{a+b}{2}\rfloor-t+1$.

{\bf Subcase 2.1.} $C_{b}$ and $C_{c}$ are odd cycles. Then $d(G-w_{0}w_{1})=\lfloor \frac {a}{2}\rfloor+b-t\geq d+1$, $d(G-w_{0}-w_{1})\geq c-3\geq d$ and $d\geq3$. By Lemmas 2.2, 2.8, 2.9, 2.11 and 2.13,
\begin{align*}
 a_{2i}(G^{\sigma})=&a_{2i}(G^{\sigma}-w_{0}w_{1})+a_{2i-2}(G^{\sigma}-w_{0}-w_{1})\\
                      \geq&a_{2i}(U^{-}_{n,d+1})+a_{2i-2}(T_{n-2,d}).
\end{align*}
Combining with $(3)$ and Lemma 2.14, then $G^{\sigma}\succ B^{-,-,-}_{n,d+1}\succ B^{-,-,-}_{n,d}$.

{\bf Subcase 2.2.} $C_{b}$ is an odd cycle and $C_{c}$ is an even cycle. If $b=5,$ we have $G^{\sigma}\succ B^{-,-,-}_{n,d}$ by Lemma 2.1. Otherwise, $d(G-w_{0}w_{1}-w_{1}w_{2})=\lfloor \frac {a}{2}\rfloor+b-t-1\geq d+1$, $d(G-w_{1}-w_{2}-w_{3}w_{4})\geq d-1$, $d(G-w_{0}-w_{1}-w_{2}w_{3})\geq c-4\geq d$ and $d\geq4$. By Lemmas 2.2, 2.8, 2.9, 2.11, 2.12 and 2.13, we have
\begin{align*}
 a_{2i}(G^{\sigma})\geq&a_{2i}(G^{\sigma}-w_{0}w_{1}-w_{1}w_{2})+a_{2i-2}(G^{\sigma}-w_{0}-w_{1})\\
                     &+a_{2i-2}(G^{\sigma}-w_{1}-w_{2})-2a_{2i-c}(G^{\sigma}-V(C_{c}))\\
                     \geq&a_{2i}(G^{\sigma}-w_{0}w_{1}-w_{1}w_{2})+a_{2i-2}(G^{\sigma}-w_{0}-w_{1}-w_{2}w_{3})\\
                     &+a_{2i-2}(G^{\sigma}-w_{1}-w_{2}-w_{3}w_{4})\\
                     \geq&a_{2i}(U^{-}_{n-1,d+1})+a_{2i-2}(U^{-}_{n-3,d-1})+a_{2i-2}(T_{n-3,d})\\
                      \geq&a_{2i}(U^{-}_{n-1,d+1})+a_{2i-2}(T_{n-3,d-1})+a_{2i-2}(T_{n-3,d}),
\end{align*}
which, together with (4) and Lemma 2.14, implies $G^{\sigma}\succ B^{-,-,-}_{n,d+1}\succ B^{-,-,-}_{n,d}$.

For $C_{c}$ is an odd cycle and $C_{b}$ is an even cycle, we have $G^{\sigma}\succ B^{-,-,-}_{n,d}$ by similar arguments as above.

{\bf Subcase 2.3.} $C_{b}$ and $C_{c}$ are even cycles. Then $c\geq b\geq 6$ and $C_{a}$ is an even cycle. $d(G-u_{0})=c-2\geq d+1$, $d\geq3$ and $n=a+b-t$. When $a\neq4$ and $b\neq6$, we have $n-d\geq5$. By Lemmas 2.3, 2.6, 2.7, 2.8 and 2.9, we have
\begin{align*}
 a_{2i}(G^{\sigma})\geq&a_{2i}(G^{\sigma}-u_{0})+a_{2i-2}(G^{\sigma}-u_{0}-u_{1})+a_{2i-2}(G^{\sigma}-u_{0}-u_{b-1})\\
                     &+a_{2i-2}(G^{\sigma}-u_{0}-v_{a-1})-2a_{2i-a}(G^{\sigma}-V(C_{a}))\\
                     &-2a_{2i-b}(G^{\sigma}-V(C_{b}))-2a_{2i-c}(G^{\sigma}-V(C_{c}))\\
                     \geq&a_{2i}(G^{\sigma}-u_{0})+a_{2i-2}(G^{\sigma}-u_{0}-u_{1}-u_{t-1})\\
                     &+a_{2i-2}(G^{\sigma}-u_{0}-u_{b-1}-u_{t-1})+a_{2i-2}(G^{\sigma}-u_{0}-v_{a-1}-u_{t-1})\\
                     \geq&a_{2i}(T_{n-1,d+1})+a_{2i-2}(P_{a-t}\cup P_{b-t-1}\cup P_{t-2})+a_{2i-2}(P_{a-t}\cup P_{b-t}\cup P_{t-3})\\
                     &+a_{2i-2}(P_{a-t-1}\cup P_{b-t}\cup P_{t-2})\quad (By\ Lemma\ 2.6)\\
                      \geq&a_{2i}(T_{n-1,d+1})+3a_{2i-2}(P_{a+b-t-5})\quad (By\ Lemma\ 2.7)\\
                      \geq&a_{2i}(T_{n-1,d+1})+a_{2i-2}(P_{n-5})+2a_{2i-2}(T_{n-5,d-1}).
\end{align*}
Combining with (5) and Lemma 2.14, $G^{\sigma}\succ B^{-,-,-}_{n,d+1}\succ B^{-,-,-}_{n,d}$. When $a=4$ and $b=6$, we have $G^{\sigma}\succ B^{-,-,-}_{n,d}$ by Lemma 2.1.    \hfill$\blacksquare$

\noindent\begin{lemma}\label{le:3.3}
Let $G\in \mathcal{B}(n,d)$ with $3\leq d\leq n-4$. If $G$ contains exactly one pendent vertex $u$ on all diametrical paths of $G$ such that $G-u$ contains no pendent vertices, then $G^{\sigma}\succ B^{-,-,-}_{n,d}$.
\end{lemma}

{\bf Proof.} Let $b\geq a$. Since $d\leq n-4$, we have $b\geq 5$. Let $v$ be the neighbor of $u$.

{\bf Case 1.} When $t\leq 1$, there are exactly two cycles $C_{a}$ and $C_{b}$ in $G$. Then $d=\lfloor \frac {a}{2}\rfloor+\lfloor \frac {b}{2}\rfloor+l(G)+1$.

{\bf Subcase 1.1.} The length of $C_{b}$ is odd. If $b\geq 7$, then $d(G-u_{1}u_{2})\geq\lfloor \frac {a}{2}\rfloor+b+l(G)-2\geq d+1$. If $b=5$ and $v$ lies on $C_{a}$, then $d(G-u_{1}u_{2})\geq\lfloor \frac {a}{2}\rfloor+b+l(G)-1=d+1$.
In these cases, by similar arguments as those in Subcase 1.1 of Lemma 3.2, $G^{\sigma}\succ B^{-,-,-}_{n,d}$. Otherwise, $a=3,4,$ $b=5$, and $v$ lies on $C_{b}$. If $t=1$, we have $G^{\sigma}\succ B^{-,-,-}_{n,d}$ by Lemma 2.1. If $t=0$, then $G$ is isomorphic to $G_{5}$ or $G_{6}$ in Figure 5.

\begin{center}
\includegraphics [width=12.5 cm, height=2 cm]{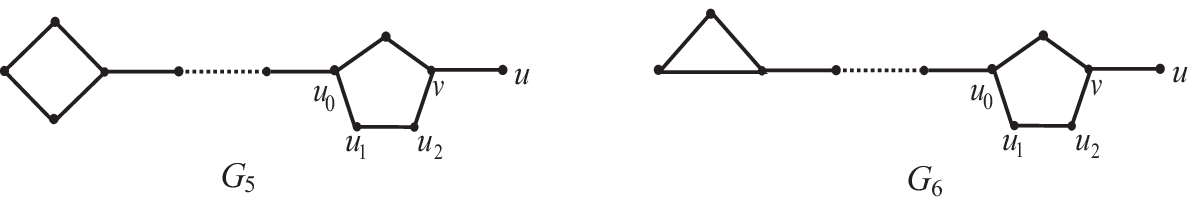}
\centerline{Figure 5: Graphs $G_{5}$ and $G_{6}$.}
\end{center}

By Lemmas 2.2, 2.11, 2.12 and 2.13, we have
\begin{align*}
 a_{2i}(G^{\sigma})=&a_{2i}(G^{\sigma}-u_{1}u_{2})+a_{2i-2}(G^{\sigma}-u_{1}-u_{2})\\
                      \geq&a_{2i}(U^{-}_{n,d})+a_{2i-2}(T_{n-2,d}),
\end{align*}
which, together with (3), implies $G^{\sigma}\succ B^{-,-,-}_{n,d}$.

{\bf Subcase 1.2.} The length of $C_{b}$ is even. If $b\geq 8$, hence $d(G-u_{1}u_{2}-u_{2}u_{3})=\lfloor \frac {a}{2}\rfloor+b+l(G)-3\geq d$. If $b=6$ and $v$ lies on $C_{a}$, then $d(G-u_{1}u_{2}-u_{2}u_{3})=\lfloor \frac {a}{2}\rfloor+b+l(G)-3=d$. In these cases, by similar arguments as those in Subcase 1.2 of Lemma 3.2, $G^{\sigma}\succ B^{-,-,-}_{n,d}$. If $a=5,$ $b=6$, and $v$ lies on $C_{b}$, then $G^{\sigma}\succ B^{-,-,-}_{n,d}$ by similar arguments as those in Subcase 1.1 of Lemma 3.2. Otherwise, $a=3,4,$ $b=6$, and $v$ lies on $C_{b}$. If $l(G)\leq2$, by Lemma 2.1, then we have $G^{\sigma}\succ B^{-,-,-}_{n,d}$ . If $l(G)\geq3$, by Lemmas 2.2, 2.10, 2.11 and 2.12, we have
\begin{align*}
 a_{2i}(G^{\sigma})\geq&a_{2i}(G^{\sigma}-u_{0}u_{1}-u_{1}u_{2})+a_{2i-2}(G^{\sigma}-u_{0}-u_{1})\\
                     &+a_{2i-2}(G^{\sigma}-u_{1}-u_{2})-2a_{2i-b}(G^{\sigma}-V(C_{b}))\\
                     \geq&a_{2i}(G^{\sigma}-u_{0}u_{1}-u_{1}u_{2})+a_{2i-2}(G^{\sigma}-u_{0}-u_{1}-u_{2}u_{3})\\
                     &+a_{2i-2}(G^{\sigma}-u_{1}-u_{2}-u_{3}u_{4})\\
                     \geq&a_{2i}(U^{-}_{n-1,d})+a_{2i-2}(U^{-}_{n-4,d-2}\cup P_{2})+a_{2i-2}(U^{-}_{n-7,d-5}\cup P_{4})\\
                      \geq&a_{2i}(U^{-}_{n-1,d})+a_{2i-2}(T_{n-4,d-2}\cup P_{2})+a_{2i-2}(T_{n-7,d-5}\cup P_{4})\\
                      \geq&a_{2i}(U^{-}_{n-1,d})+a_{2i-2}(T_{n-3,d-1})+a_{2i-2}(T_{n-4,d-2}).
\end{align*}
Combining with (4), $G^{\sigma}\succ B^{-,-,-}_{n,d}$.

{\bf Case 2.} $t\geq 2$. Then $d=\lfloor \frac{c}{2}\rfloor+1=\lfloor \frac{a+b}{2}\rfloor-t+2$. Since $b\geq 5$, assume that $w_{0}, w_{1}\neq v$. Note that $a-t+1\geq t-1$ and $b-t+1\geq t-1$. By similar arguments as those in Case 2 of Lemma 3.2, we have $G^{\sigma}\succ B^{-,-,-}_{n,d}$. \hfill$\blacksquare$


\noindent\begin{theorem}\label{th:3.4}
Let $G\in \mathcal{B}(n,d)$ with $3\leq d\leq n-3$ and $G\neq B_{n,d}$. If $t=0$, then $G^{\sigma}\succ B^{-,-,-}_{n,d}$.
\end{theorem}

{\bf Proof.} We prove this theorem by induction on $n-d$.

By Lemma 3.1, the result holds for $n-d=3$. Let $h\geq 4$ and suppose that the
result holds for $n-d<h$. Now suppose that $n-d=h$ and $G\in \mathcal{B}(n,d)$.

{\bf Case 1.} There is no pendent vertex in $G$. By Lemma 3.2, we have $G^{\sigma}\succ B^{-,-,-}_{n,d}$.

{\bf Case 2.} There is a pendent vertex $u$ outside some diametrical path $P(G)=x_{0}x_{1}\cdots x_{d}$. Let $v$ be the neighbor of $u$. Then $G-u\in \mathcal{B}(n-1,d)$. By the induction hypothesis, $G^{\sigma}-u\succ B^{-,-,-}_{n-1,d}$. By Lemma 2.2, we have
$$a_{2i}(B^{-,-,-}_{n,d})=a_{2i}(B^{-,-,-}_{n-1,d})+a_{2i-2}(T_{d+1,d-2}).\eqno{(7)}$$
Let $H=G-u-v$, it suffices to prove that $a_{2i}(H^{\sigma})\geq a_{2i}(T_{d+1,d-2})$.

{\bf Subcase 2.1.} $v$ lies on some cycle, say $C_{a}$.

{\bf Subcase 2.1.1.} Suppose that $P(G)$ and $C_{b}$ have no common vertices. Then $H\supseteq P_{k}\cup P_{d-k}\cup C_{b}$. By Lemmas 2.4, 2.5, 2.6 and 2.7, we have
\begin{align*}
 a_{2i}(H^{\sigma})&\geq a_{2i}(P_{k}\cup P_{d-k}\cup C_{b})\geq a_{2i}(P_{d-1}\cup S_{b})\geq a_{2i}(P_{d-1}\cup P_{3})\\
                      &\geq a_{2i}(P_{d+1})\geq a_{2i}(T_{d+1,d-2}).
\end{align*}

{\bf Subcase 2.1.2.} Suppose that $P(G)$ and $C_{b}$ have common vertices $x_{l},\cdots ,x_{l+q}$, where $q\geq 0$.

If $v$ lies outside $P(G)$, then $H\supseteq H_{1}$, where $H_{1}\in U(s_{1},d), s_{1}\geq d+2$. By Lemmas 2.4, 2.7, 2.11 and 2.12, we have
$$a_{2i}(H^{\sigma})\geq a_{2i}(H_{1}^{\sigma})\geq a_{2i}(U^{\sigma}_{s_{1},d})\geq a_{2i}(T_{s_{1},d})\geq a_{2i}(P_{d+1})\geq a_{2i}(T_{d+1,d-2}).$$

If $v$ lies on $P(G)$. Then $P(G)$ and $C_{a}$ have common vertices $x_{k},\cdots, x_{k+p}$, where $p\geq 0, k+p<l$.

Suppose that $p=0$, then $k\geq 1$, $H\supseteq P_{2}\cup P_{k}\cup H_{2}$, where $H_{2}\in U(s_{2},d_{2}), s_{2}\geq d_{2}+2$, $d_{2}\geq d-k-1\geq1$. If $d_{2}=1$, then $k=d-2$ and $H_{2}=C_{3}$. By Lemmas 2.4, 2.5, 2.6 and 2.7, we have
$$a_{2i}(H^{\sigma})\geq a_{2i}(P_{2}\cup P_{k}\cup C_{3})\geq a_{2i}(P_{d-1}\cup P_{3})\geq a_{2i}(P_{d+1})\geq a_{2i}(T_{d+1,d-2}).$$
If $d_{2}=2$, then $k\geq d-3$ and $s_{2}\geq 4$. By Lemmas 2.4, 2.5, 2.6, 2.8 and 2.9, we have
\begin{align*}
 a_{2i}(H^{\sigma})&\geq a_{2i}(P_{2}\cup P_{k}\cup H_{2}^{\sigma})\geq a_{2i}(P_{k+1}\cup S_{s_{2}})\geq a_{2i}(P_{k+1}\cup S_{4})\\
                      &\geq a_{2i}(T_{k+4,k+2})\geq a_{2i}(T_{d+1,d-1})\geq a_{2i}(T_{d+1,d-2}).
\end{align*}
If $d_{2}\geq 3$, then by Lemmas 2.4, 2.6, 2.8, 2.9, 2.11 and 2.12, we have
\begin{align*}
 a_{2i}(H^{\sigma})&\geq a_{2i}(P_{2}\cup P_{k}\cup H_{2}^{\sigma})\geq a_{2i}(P_{k+1}\cup U^{-}_{s_{2},d_{2}})\geq a_{2i}(P_{k+1}\cup T_{s_{2},d_{2}})\\
                      &\geq a_{2i}(T_{s_{2}+k,d_{2}+k})\geq a_{2i}(T_{d_{2}+k+2,d_{2}+k})\geq a_{2i}(T_{d+1,d-1})\geq a_{2i}(T_{d+1,d-2}).
\end{align*}





Suppose that $p\geq1$. If $v=x_{k},$ then $H\supseteq P_{k}\cup H_{3}$, where $H_{3}\in U(s_{3},d_{3}), s_{3}\geq d_{3}+2$, $d_{3}\geq d-k\geq3$. By Lemmas 2.4, 2.10, 2.11 and 2.12, we have
$$a_{2i}(H^{\sigma})\geq a_{2i}(P_{k}\cup U_{s_{3},d_{3}}^{-})\geq a_{2i}(P_{k}\cup T_{s_{3},d_{3}})\geq a_{2i}(T_{s_{3}+k-1,d_{3}+k-1})\geq a_{2i}(T_{d+1,d-2}).$$

If $v=x_{k+p}$, $H\supseteq P_{k+p+1}\cup H_{4}$ or $T_{1}\cup H_{4}$($p\geq2$), where $H_{4}\in U(s_{4},d_{4}), s_{4}\geq d_{4}+2$, $d_{4}\geq d-k-p-1\geq1$ and $T_{1}\in T(k+p+1,k+p-1)$. If $d_{4}=1$, then $k+p=d-2$ and $H_{4}=C_{3}$. By Lemmas 2.4, 2.5, 2.6, 2.7, 2.8, 2.9 and 2.10, we have
$$a_{2i}(H^{\sigma})\geq a_{2i}(P_{d-1}\cup C_{3})\geq a_{2i}(P_{d-1}\cup P_{2})\geq a_{2i}(P_{d+1})\geq a_{2i}(T_{d+1,d-2}).$$
or
$$a_{2i}(H^{\sigma})\geq a_{2i}(T_{d-1,d-3}\cup C_{3})\geq a_{2i}(T_{d-1,d-3}\cup P_{3})\geq a_{2i}(T_{d+1,d-2}).$$
If $d_{4}=2$, then $k+p\geq d-3$ and $s_{4}\geq 4$. By Lemmas 2.4, 2.5, 2.8, 2.9 and 2.10, we have
$$a_{2i}(H^{\sigma})\geq a_{2i}(P_{k+p+1}\cup S_{s_{4}})\geq a_{2i}(P_{k+p+1}\cup T_{4,2})\geq a_{2i}(T_{d+1,d-2}).$$
or
$$a_{2i}(H^{\sigma})\geq a_{2i}(T_{k+p+1,k+p-1}\cup S_{s_{4}})\geq a_{2i}(T_{k+p+1,k+p-1}\cup T_{4,2})\geq a_{2i}(T_{d+1,d-2}).$$
If $d_{4}\geq 3$, then by Lemmas 2.4, 2.8, 2.9, 2.10, 2.11 and 2.12, we have
$$a_{2i}(H^{\sigma})\geq a_{2i}(P_{k+p+1}\cup U^{-}_{s_{4},d_{4}})\geq a_{2i}(P_{k+p+1}\cup T_{s_{4},d_{4}})\geq a_{2i}(T_{d+1,d-2}).$$
or
$$a_{2i}(H^{\sigma})\geq a_{2i}(T_{k+p+1,k+p-1}\cup U^{-}_{s_{4},d_{4}})\geq a_{2i}(T_{k+p+1,k+p-1}\cup T_{s_{4},d_{4}})\geq a_{2i}(T_{d+1,d-2}).$$

{\bf Subcase 2.2.} $v$ lies outside any cycle.

{\bf Subcase 2.2.1.} Suppose that $v$ lies on $P(G)$ and $v=x_{k}$.

If $P(G)$ and any cycle have no common vertices, then $H\supseteq C_{a}\cup C_{b}\cup P_{k}\cup P_{d-k}$. By Lemmas 2.4, 2.5, 2.6 and 2.7, we have
$$a_{2i}(H^{\sigma})\geq a_{2i}(C_{a}\cup C_{b}\cup P_{k}\cup P_{d-k})\geq a_{2i}(P_{3}\cup P_{d-1}))\geq a_{2i}(T_{d+1,d-2}).$$

If $P(G)$ and exactly one cycle have no common vertices, say $C_{a}$, then $H\supseteq C_{a}\cup P_{k}\cup H_{1}$, where $H_{1}\in U(s_{1},d_{1}), s_{1}\geq d_{1}+2, d_{1}\geq d-k-1\geq1$. If $d_{1}=1$, then $k=d-2$ and $H_{1}=C_{3}$. By Lemmas 2.4, 2.5, 2.6 and 2.7, we have
$$a_{2i}(H^{\sigma})\geq a_{2i}(P_{2}\cup P_{k}\cup C_{3})\geq a_{2i}(P_{d-1}\cup P_{2})\geq a_{2i}(P_{d+1})\geq a_{2i}(T_{d+1,d-2}).$$
If $d_{1}\geq2$, then $k\geq d-3$. By Lemmas 2.4, 2.5, 2.6 and 2.10, we have
$$a_{2i}(H^{\sigma})\geq a_{2i}(C_{a}\cup P_{k}\cup H_{1}^{\sigma})\geq a_{2i}( P_{k+1}\cup S_{s_{1}})\geq a_{2i}(T_{d+1,d-2}).$$

If $P(G)$ and two cycles have common vertices, then $H\supseteq P_{k}\cup H_{2}$ or $H_{3}\cup H_{4}$, where $H_{2}\in B(s_{2},d_{2})$, $H_{3}\in U(s_{3},d_{3})$, $H_{4}\in U(s_{4},d_{4})$, $d_{2}+3\leq s_{2}\leq n-k-2$, $d_{2}\geq d-k-1\geq4$, $s_{3}\geq d_{3}+2$, $d_{3}\geq k-1\geq2$, $s_{4}\geq d_{4}+2$, $d_{4}\geq d-k-1\geq1$.

Suppose that $H\supseteq P_{k}\cup H_{2}$, $s_{2}-d_{2}<h$ and $d_{2}\geq4$, by the induction hypothesis, $H_{2}^{\sigma}\succ B_{s_{2},d_{2}}$. By Lemmas 2.4 and 2.10, and (6), we have
$$a_{2i}(H^{\sigma})\geq a_{2i}(P_{k}\cup H^{\sigma}_{2})\geq a_{2i}(P_{k}\cup B^{-,-,-}_{s_{2},d_{2}})\geq a_{2i}(P_{k}\cup T_{s_{2},d_{2}})\geq a_{2i}(T_{d+1,d-2}).$$

Suppose that $H\supseteq H_{3}\cup H_{4}$.  If $d_{3}=2$ and $d_{4}=1$, then $d=5$ and $s_{3}\geq4$. By Lemmas 2.4, 2.5, 2.9 and 2.10, we have
$$a_{2i}(H^{\sigma})\geq a_{2i}(H^{\sigma}_{3}\cup H^{\sigma}_{4})\geq a_{2i}(S_{s_{3}}\cup C^{-}_{3})\geq a_{2i}(P_{3}\cup T_{4,2})\geq a_{2i}(T_{6,4})\geq a_{2i}(T_{d+1,d-2}).$$
If $d_{3}=2$ and $d_{4}=2$, then $d=6$ and $s_{3}\geq4, s_{4}\geq4$. By Lemmas 2.4, 2.5 and 2.10, we have
$$a_{2i}(H^{\sigma})\geq a_{2i}(H^{\sigma}_{3}\cup H^{\sigma}_{4})\geq a_{2i}(S_{s_{3}}\cup S_{s_{4}})\geq a_{2i}(T_{4,2}\cup T_{4,2})\geq a_{2i}(T_{d+1,d-2}).$$
If $d_{3}\geq3$ and $d_{4}=1$, then $d_{3}\geq d-3$. By Lemmas 2.4, 2.5, 2.9, 2.10, 2.11 and 2.12, we have
$$a_{2i}(H^{\sigma})\geq a_{2i}(H^{\sigma}_{3}\cup H^{\sigma}_{4})\geq a_{2i}(U^{-}_{s_{3},d_{3}}\cup C^{-}_{3})\geq a_{2i}(T_{s_{3},d_{3}}\cup P_{3})\geq a_{2i}(T_{d+1,d-2}).$$
If $d_{3}\geq3$ and $d_{4}=2$, then $d_{3}\geq d-4, s_{4}\geq4$. By Lemmas 2.4, 2.5, 2.10, 2.11 and 2.12, we have
$$a_{2i}(H^{\sigma})\geq a_{2i}(H^{\sigma}_{3}\cup H^{\sigma}_{4})\geq a_{2i}(U^{-}_{s_{3},d_{3}}\cup S_{s_{4}})\geq a_{2i}(T_{s_{3},d_{3}}\cup T_{4,2})\geq a_{2i}(T_{d+1,d-2}).$$
If $d_{3}\geq3$ and $d_{4}\geq3$, then $d_{3}+d_{4}\geq d-2$. By Lemmas 2.4, 2.10, 2.11 and 2.12, we have
$$a_{2i}(H^{\sigma})\geq a_{2i}(H^{\sigma}_{3}\cup H^{\sigma}_{4})\geq a_{2i}(U^{-}_{s_{3},d_{3}}\cup U^{-}_{s_{4},d_{4}})\geq a_{2i}(T_{s_{3},d_{3}}\cup T_{s_{4},d_{4}})\geq a_{2i}(T_{d+1,d-2}).$$

{\bf Subcase 2.2.2.} Suppose that $v$ lies outside $P(G)$. Then $G\supseteq C_{a}\cup C_{b}\cup P(G)$, $C_{a}\cup H_{1}$ or $H_{2}$, where $H_{1}\in U(s_{1},d)$ with $s_{1}\geq d+2$ and $H_{2}\in B(s_{2},d)$ with $d+3\leq s_{2}\leq n-2$. We can prove $a_{2i}(H^{\sigma})\geq a_{2i}(T_{d+1,d-2})$ by similar arguments as above.

{\bf Case 3.} All pendent vertices are contained in the $P(G)$, where $P(G)=x_{0}x_{1}\cdots x_{d}$ is a diametrical path of $G$. Suppose that $y_{0}y_{1}\cdots y_{p}$ is a path whose internal vertices $y_{1},y_{2},\cdots, y_{p-1}$ all have degree two and $y_{p}$ is a pendent vertex. Then
we say that it is a pendent path, denoted by $(y_{0},y_{p})$ (see \cite{YZY}).

{\bf Subcase 3.1.} There are exactly two pendent vertices, $x_{0}$ and $x_{d}$. Suppose that $deg(x_{k}),$ $deg(x_{l})\geq3$ and that $(x_{k},x_{0})$ and $(x_{l},x_{d})$ are distinct pendent paths. Let $s=k-l$.

If $s=0 (x_{k}=x_{l})$, then $k\geq3$ and $l<d-3$. So it suffices
to prove that $H^{\sigma}_{1}, H^{\sigma}_{2}\succ B^{-,-,-}_{n-d+3,3}$, $H^{\sigma}_{3}, H^{\sigma}_{4}\succ S_{n-d+2}$, where $H_{1}=G-\{x_{k-3},\cdots , x_{0}\}-\{x_{l+2},\cdots,x_{d}\}$, $H_{3}=G-\{x_{k-3},\cdots,x_{0}\}-\{x_{l+1},\cdots,x_{d}\}$, $H_{2}=G-\{x_{k-2},\cdots , x_{0}\}-\{x_{l+3},\cdots,x_{d}\}$, and $H_{4}=G-\{x_{k-2},\cdots , x_{0}\}-\{x_{l+2},\cdots,x_{d}\}$. By Lemma 2.5, $H^{\sigma}_{3}, H^{\sigma}_{4}\succ S_{n-d+2}$. Let $d_{1}=d(H_{1})$. Since $d_{1}\geq4, n-d+3-d_{1}<h$. $H^{\sigma}_{1}\succ B^{-,-,-}_{n-d+3,d_{1}}\succ B^{-,-,-}_{n-d+3,3}$ by the induction hypothesis and Lemma 2.14. Similarly, $H^{\sigma}_{2}\succ B^{-,-,-}_{n-d+3,3}$.

If $s=1$ or 2, then by similar arguments as above, we have $G^{\sigma}\succ B^{-,-,-}_{n,d}$.

If $s\geq3$, then we only need consider the case $k\geq2$ and $l\leq d-2$. So it suffices
to prove that $H^{\sigma}_{5}, H^{\sigma}_{6}\succ B^{-,-,-}_{n-d+s+1,s+1}$, $H^{\sigma}_{7}\succ B^{-,-,-}_{n-d+s+2,s+2}$, $H^{\sigma}_{8}\succ B^{-,-,-}_{n-d+s,s}$, where $H_{5}=G-\{x_{k-2},\cdots,x_{0}\}-\{x_{l+1},\cdots,x_{d}\}$, $H_{7}=G-\{x_{k-2},\cdots,x_{0}\}-\{x_{l+2},\cdots,x_{d}\}$, $H_{8}=G-\{x_{k-1},\cdots,x_{0}\}-\{x_{l+1},\cdots,x_{d}\}$, and $H_{6}=G-\{x_{k-1},\cdots,x_{0}\}-\{x_{l+2},\cdots,x_{d}\}$. Let $d_{j}=d(H_{j})$ and $n_{j}=|V(H_{j})|$, where $j=5,6,7,8$. Then $d_{j}\geq 4$. If $n_{j}-d_{j}<h$, then by the induction hypothesis and Lemma 2.14, we have the desired result.

Suppose that $n_{j}-d_{j}=h$. If $x_{k-1}$ lies on all diametrical paths of $H_{5}$, then by Lemmas 2.14 and 3.3, $H^{\sigma}_{5}\succ B^{-,-,-}_{n-d+s+1,s+1}$. Otherwise, by similar arguments as those in Case 2, we also have $H^{\sigma}_{5}\succ B^{-,-,-}_{n-d+s+1,s+1}$. Similarly, $H^{\sigma}_{6}\succ B^{-,-,-}_{n-d+s+1,s+1}$. By Lemmas 2.14 and 3.2, we have $H^{\sigma}_{8}\succ B^{-,-,-}_{n-d+s,s}$. If there exists some diametrical path $P(H_{7})$
such that $x_{k-1}$ or $x_{l+1}$ lies outside $P(H_{7})$, then by similar arguments as those in Case 2, we have $H^{\sigma}_{7}\succ B^{-,-,-}_{n-d+s+2,s+2}$. Otherwise, by Lemmas 2.2, 2.4, 2.11, 2.12, 2.13 and 3.3, we have
$H^{\sigma}_{7}-x_{k-1}\succ B^{-,-,-}_{n-d+s+1,s+2}$, $H^{\sigma}_{7}-x_{k-1}-x_{k}\succ U^{-}_{n-d+s,s}\succ T_{n-d+s,s}\succ T_{s+3,s}$, and then $H^{\sigma}_{7}\succ B^{-,-,-}_{n-d+s+2,s+2}$.

{\bf Subcase 3.2.} There is only one pendent vertex. By similar arguments as those in Subcase 3.1, we have $G^{\sigma}\succ B^{-,-,-}_{n,d}$.

Combining all those cases above, we complete the proof.     \hfill$\blacksquare$

\noindent\begin{theorem}\label{th:3.5}
Let $G\in \mathcal{B}(n,d)$ with $3\leq d\leq n-3$ and $G\neq B_{n,d}$. If $t\geq1$, then $G^{\sigma}\succ B^{-,-,-}_{n,d}$.
\end{theorem}

{\bf Proof.} We prove this theorem by induction on $n-d$.

By Lemma 3.1, the result holds for $n-d=3$. Let $h\geq 4$ and suppose that the
result holds for $n-d<h$. Now suppose that $n-d=h$ and $G\in \mathcal{B}(n,d)$.

{\bf Case 1.} There is no pendent vertex in $G$. By Lemma 3.2, we have $G^{\sigma}\succ B^{-,-,-}_{n,d}$.

{\bf Case 2.} There is a pendent vertex $u$ outside some diametrical path $P(G)=x_{0}x_{1}\cdots x_{d}$. Let $v$ be the neighbor of $u$. Then $G-u\in \mathcal{B}(n-1,d)$. By the induction hypothesis, $G^{\sigma}-u\succ B^{-,-,-}_{n-1,d}$. Let $H=G-u-v$, it suffices to prove that $a_{2i}(H^{\sigma})\geq a_{2i}(T_{d+1,d-2})$ by the $(7)$.

{\bf Subcase 2.1.} $v$ lies on some cycle, say $C_{a}$.

{\bf Subcase 2.1.1.} Suppose that $v=u_{0}$ or $u_{t-1}$.

If $v$ lies outside $P(G)$, then $H\supseteq P(G)$. By Lemmas 2.4 and 2.7, we have $a_{2i}(H^{\sigma})\geq a_{2i}(P_{d+1})\geq a_{2i}(T_{d+1,d-2})$

If $v$ lies on $P(G)$. Let $v=x_{k}$. If $C_{a}$ and $C_{b}$ have exactly one common vertex, then $H\supseteq P_{2}\cup P_{2}\cup P_{k}\cup P_{d-k}$, $P_{2}\cup P_{k}\cup P_{d-k+1}$, $P_{2}\cup P_{k}\cup T_{1}$, $P_{k+1}\cup P_{d-k+1}$, $P_{k+1}\cup T_{1}$ or $T_{1}\cup T_{2}$, where $T_{1}\in T(d-k+1,d-k-1)$, $T_{2}\in T(k+1,k-1)$. If $C_{a}$ and $C_{b}$ have at least two common vertices, $H\supseteq P_{3}\cup P_{k}\cup P_{d-k}$, $P_{k}\cup P_{d-k+2}$, $P_{k}\cup T_{3}$, $P_{k}\cup T_{4}$ or $P(G)$, where $T_{3}\in T(d-k+2,d-k-1)$, $T_{4}\in T(d-k+2,d-k)$.

If $H\supseteq P_{2}\cup P_{2}\cup P_{k}\cup P_{d-k}$, $P_{2}\cup P_{k}\cup P_{d-k+1}$, $P_{k+1}\cup P_{d-k+1}$, $P_{3}\cup P_{k}\cup P_{d-k}$, $P_{k}\cup P_{d-k+2}$ or $P(G)$, by Lemmas 2.4, 2.6 and 2.7, then we have
$$a_{2i}(H^{\sigma})\geq a_{2i}(P_{d+1})\geq a_{2i}(T_{d+1,d-2}).$$
If $H\supseteq P_{2}\cup P_{k}\cup T_{1}$ or $P_{k+1}\cup T_{1}$, by Lemmas 2.4, 2.6, 2.8, 2.9 and 2.10, then we have
$$a_{2i}(H^{\sigma})\geq a_{2i}(P_{k+1}\cup T_{d-k+1,d-k-1})\geq a_{2i}(T_{d+1,d-1})\geq a_{2i}(T_{d+1,d-2}).$$
If $H\supseteq T_{1}\cup T_{2}$, by Lemmas 2.4, 2.8 and 2.10, then we have
$$a_{2i}(H^{\sigma})\geq a_{2i}(T_{k+1,k-1}\cup T_{d-k+1,d-k-1})\geq a_{2i}(T_{d+1,d-2}).$$
If $H\supseteq P_{k}\cup T_{3}$, by Lemmas 2.4, 2.8 and 2.10, then we have
$$a_{2i}(H^{\sigma})\geq a_{2i}(P_{k}\cup T_{d-k+2,d-k-1})\geq a_{2i}(T_{d+1,d-2}).$$
If $H\supseteq P_{k}\cup T_{4}$, by Lemmas 2.4, 2.8, 2.9 and 2.10, then we have
$$a_{2i}(H^{\sigma})\geq a_{2i}(P_{k}\cup T_{d-k+2,d-k})\geq a_{2i}(T_{d+1,d-1})\geq a_{2i}(T_{d+1,d-2}).$$

{\bf Subcase 2.1.2.} Suppose that $v\neq u_{0}$ and $u_{t-1}$. If $v$ lies outside $P(G)$, then $H\supseteq H_{1}$ or $P(G)\cup C_{s}$, where $H_{1}\in U(s_{1},d)$, $s_{1}\geq d+2$ $s=b$ or $c$.

If $H\supseteq H_{1}$, by Lemmas 2.4, 2.7, 2.11 and 2.12, then we have
$$a_{2i}(H^{\sigma})\geq a_{2i}(H^{\sigma}_{1})\geq a_{2i}(U^{-}_{s_{1},d})\geq a_{2i}(T_{s_{1},d})\geq a_{2i}(P_{d+1})\geq a_{2i}(T_{d+1,d-2}).$$
If $H\supseteq P(G)\cup C_{s}$, by Lemmas 2.4, 2.5 and 2.7, then we have
$$a_{2i}(H^{\sigma})\geq a_{2i}(P_{d+1}\cup S_{s})\geq a_{2i}(P_{d+1})\geq a_{2i}(T_{d+1,d-2}).$$

If $v$ lies on $P(G)$, then $P(G)$ and $C_{a}$ have common vertices, say, $x_{k},\cdots,x_{k+p}$, where $p\geq0$.

If $p=0$, then $k\geq1$, $H\supseteq P_{k}\cup P_{d-k}\cup C_{b}$. By Lemmas 2.4, 2.5, 2.6 and 2.7, then we have
$$a_{2i}(H^{\sigma})\geq a_{2i}(P_{k}\cup P_{d-k}\cup S_{b})\geq a_{2i}(P_{d-1}\cup P_{3})\geq a_{2i}(P_{d+1})\geq a_{2i}(T_{d+1,d-2}).$$

If $p\geq1$. If $v\neq x_{k}, x_{k+p}$, then $H\supseteq H_{2}$, where $H_{2}\in U(s_{2},d)$, $s_{2}\geq d+2$. We can prove that $a_{2i}(H^{\sigma})\geq a_{2i}(T_{d+1,d-2})$ by the similar arguments as above. For $v= x_{k}$ or $v= x_{k+p}$, say $v= x_{k}$, then $k\geq1$, $H\supseteq P_{k}\cup H_{3}$ or $P_{k}\cup H_{4}$, where $H_{3}\in U(s_{3},d_{3})$, $s_{3}\geq d_{3}+2$, $d_{3}\geq d-k\geq2$ and $H_{4}$ is a graph obtained by attaching $P_{d-k-2}$ to a vertex of $C_{3}$.

Suppose that $H\supseteq P_{k}\cup H_{3}$. If $d_{3}=2$, then $k=d-2$, $s_{3}\geq4$. By Lemmas 2.4, 2.5, 2.9 and 2.10, we have
$$a_{2i}(H^{\sigma})\geq a_{2i}(P_{k}\cup S_{s_{3}})\geq a_{2i}(P_{d-2}\cup T_{4,2})\geq a_{2i}(T_{d+1,d-1})\geq a_{2i}(T_{d+1,d-2}).$$
If $d_{3}\geq3$, by Lemmas 2.4, 2.9, 2.10, 2.11 and 2.12, then we have
$$a_{2i}(H^{\sigma})\geq a_{2i}(P_{k}\cup U^{-}_{s_{3},d_{3}})\geq a_{2i}(P_{k}\cup T_{s_{3},d_{3}})\geq a_{2i}(T_{d+1,d-1})\geq a_{2i}(T_{d+1,d-2}).$$

Suppose that $H\supseteq P_{k}\cup H_{4}$. If $d-k-2=0$, then $k=d-2$. We also have $a_{2i}(C_{3})\geq a_{2i}(T_{4,2})$. By Lemmas 2.4, 2.9 and 2.10, we have
$$a_{2i}(H^{\sigma})\geq a_{2i}(P_{k}\cup C_{3})\geq a_{2i}(P_{d-2}\cup T_{4,2})\geq a_{2i}(T_{d+1,d-1})\geq a_{2i}(T_{d+1,d-2}).$$
If $d-k-2\geq 1$, by Lemmas 2.2, 2.4, 2.6 and 2.8, then we have
\begin{align*}
 a_{2i}(H^{\sigma})&\geq a_{2i}(P_{k}\cup H_{4}^{*}-u_{0}u_{1})+a_{2i-2}(P_{k}\cup H_{4}^{*}-u_{0}-u_{1})\\
 &\geq a_{2i}(P_{k}\cup T_{d-k+1,d-k-1})+a_{2i-2}(P_{k}\cup P_{d-k-1})\\
 &=a_{2i}(P_{k}\cup P_{d-k-1})+2a_{2i-2}(P_{k}\cup P_{d-k-2})+a_{2i-2}(P_{k}\cup P_{d-k-1})\\
 &\geq a_{2i}(P_{d-2})+3a_{2i-2}(P_{d-3})=a_{2i}(T_{d+1,d-2}).
\end{align*}

{\bf Subcase 2.2.} $v$ lies outside any cycle.

{\bf Subcase 2.2.1} Suppose that $v$ lies on $P(G)$. Let $v=x_{k}$. If $P(G)$ and any cycle have no common vertices, then $H\supseteq C_{c}\cup P_{k}\cup P_{d-k}$. By Lemmas 2.4, 2.5, 2.6 and 2.7, we have
$$a_{2i}(H^{\sigma})\geq a_{2i}(C_{c}\cup P_{k}\cup P_{d-k})\geq a_{2i}(P_{3}\cup P_{d-1}))\geq a_{2i}(T_{d+1,d-2}).$$

If some vertex of $P(G)$ lies on one cycle, then $H\supseteq P_{k}\cup H_{1}$, where $H_{1}\in B(s_{1},d_{1})$, $d_{1}+2\leq s_{1}\leq n-2-k$, $d_{1}\geq max\{d-k-1,2\}$. Suppose that $s_{1}\geq d_{1}+3$. If $d_{1}=2$, then $k\geq d-3$, $s_{1}\geq 5$. By Lemmas 2.4, 2.5 and 2.10, we have
$$a_{2i}(H^{\sigma})\geq a_{2i}(P_{k}\cup S_{s_{1}})\geq a_{2i}(P_{k}\cup T_{5,2})\geq a_{2i}(T_{d+1,d-2}).$$
If $d_{1}\geq3$, then $s_{1}-d_{1}<h$, $H^{\sigma}_{1}\succ B^{-,-,-}_{s_{1},d_{1}}$ by the induction hypothesis. So by Lemmas 2.4 and 2.10, and (6), we have
$$a_{2i}(H^{\sigma})\geq a_{2i}(P_{k}\cup H^{\sigma}_{1})\geq a_{2i}(P_{k}\cup B^{-,-,-}_{s_{1},d_{1}})\geq a_{2i}(P_{k}\cup T_{s_{1},d_{1}})\geq a_{2i}(T_{d+1,d-2}).$$

Suppose that $s_{1}=d_{1}+2$. Then $H_{1}$ is obtained by attaching respectively paths $P_{l}$ and $P_{d_{1}-l-2}$ to the two non-adjacent vertices in $K_{4}-e$. If $d_{1}=2$, then $k\geq d-3$. We also have $a_{2i}((K_{4}-e)^{*,*,-})\geq a_{2i}(T_{5,2})$, so by Lemmas 2.4 and 2.10, we have
$$a_{2i}(H^{\sigma})\geq a_{2i}(P_{k}\cup (K_{4}-e)^{*,*,-})\geq a_{2i}(P_{k}\cup T_{5,2})\geq a_{2i}(T_{d+1,d-2}).$$
If $d_{1}\geq3$, by Lemmas 2.2, 2.4, 2.11 and 2.12, then we have
\begin{align*}
 a_{2i}(H^{\sigma})&\geq a_{2i}(P_{k}\cup H_{1}^{*,*,-}-u_{0}u_{1})+a_{2i-2}(P_{k}\cup H_{1}^{*,*,-}-u_{0}-u_{1})\\
 &\geq a_{2i}(P_{k}\cup U^{-}_{s_{1},d_{1}})+a_{2i-2}(P_{k}\cup P_{l+1}\cup P_{d_{1}-l-1})\\
 &\geq a_{2i}(P_{k}\cup T_{d-k+1,d-k-1})+a_{2i-2}(P_{k}\cup P_{d-k-2})\geq a_{2i}(T_{d+1,d-2}).
\end{align*}

{\bf Subcase 2.2.2.} Suppose that $v$ lies outside $P(G)$, $G\supseteq C_{a}\cup C_{b}\cup P(G)$ or $H_{1}$, where $H_{1}\in U(s_{1},d)$  with $d+2\leq s_{2}\leq n-2$. We can prove $a_{2i}(H^{\sigma})\geq a_{2i}(T_{d+1,d-2})$ by similar arguments as above.

{\bf Case 3.} All pendent vertices are contained in the $P(G)$. By similar arguments
as those in Case 3 of Theorem 3.4, $G^{\sigma}\succ B^{-,-,-}_{n,d}$.

Combining all those cases above, we complete the proof. \hfill$\blacksquare$



\end{document}